\documentclass[11pt]{article}

\usepackage{amssymb}
\usepackage{amsmath,amsthm}
\usepackage{amsfonts}
\usepackage[dvips]{epsfig}
\usepackage[font=small,labelfont=bf]{caption}
\usepackage{color}
\usepackage{graphicx}
\usepackage{subcaption}
\usepackage{bm}

\usepackage{xcolor}

\newtheorem{Lemma}{Lemma}[section]
\newtheorem{Theorem}{Theorem}

\makeatletter\@addtoreset{figure}{section}\makeatother

\makeatletter \@addtoreset{equation}{section} \makeatother

\newcommand{\R}{\mathbb{R}}
\newcommand{\C}{\mathbb{C}}

\renewcommand{\leq}{\leqslant}
\renewcommand{\geq}{\geqslant}
\newfam\bifam
\font\tenbi=cmmib10 scaled \magstep1 \font\sevenbi=cmmib10 at 11pt
\font\fivebi=cmmib10 at 6pt \textfont\bifam = \tenbi
\scriptfont\bifam = \sevenbi \scriptscriptfont\bifam= \fivebi

\begin{document}

\title{\textbf{Heteroclinic orbits for a system of amplitude equations for orthogonal domain walls}}
\author{Boris Buffoni\\
{\footnotesize Ecole Polytechnique F\'ed\'erale de Lausanne, Institut de Math\'ematiques,}\\
{\footnotesize   Station 8, 1015 Lausanne, Switzerland}\\
\and Mariana Haragus \\
{\footnotesize FEMTO-ST institute, Univ. Bourgogne Franche-Comt\'e, CNRS,}\\
{\footnotesize 15B avenue des Montboucons, 25030 Besan\c con cedex, France}
\and G\'{e}rard Iooss \\
{\footnotesize Laboratoire J.A.Dieudonn\'e, I.U.F., Universit\'e C\^ote d'Azur, CNRS,}\\
{\footnotesize Parc Valrose, 06108 Nice cedex 2, France} }
\date{\today}
\maketitle

\begin{abstract}
Using a variational method, we prove the existence of heteroclinic solutions for a $6$-dimensional system of ordinary differential equations. We derive this system from the classical B\'enard-Rayleigh problem near the convective instability threshold. The constructed heteroclinic solutions provide first order approximations for domain walls between two orthogonal convective rolls.   
\end{abstract}

\section{Introduction}

We consider the following system of ordinary differential equations
\begin{eqnarray}
\frac{d^{4}A_{0}}{dx^{4}} &=&A_{0}(1-|A_{0}|^{2}-g|B_{0}|^{2}),
\label{e:A0} \\
\frac{d^{2}B_{0}}{dx^{2}} &=&\varepsilon
^{2}B_{0}(-1+g|A_{0}|^{2}+|B_{0}|^{2}),  \label{e:B0}
\end{eqnarray}
in which $A_0$ and $B_0$ are complex-valued functions defined on $\R$, $\varepsilon$ is a nonnegative parameter and $g>1$. The purpose of this work is twofold: to rigorously derive this system from the B\'enard-Rayleigh convection problem and to prove that it possesses heteroclinic orbits.

\begin{figure}[th]
  \centering
  \vspace*{2ex}
\includegraphics[width=16ex]{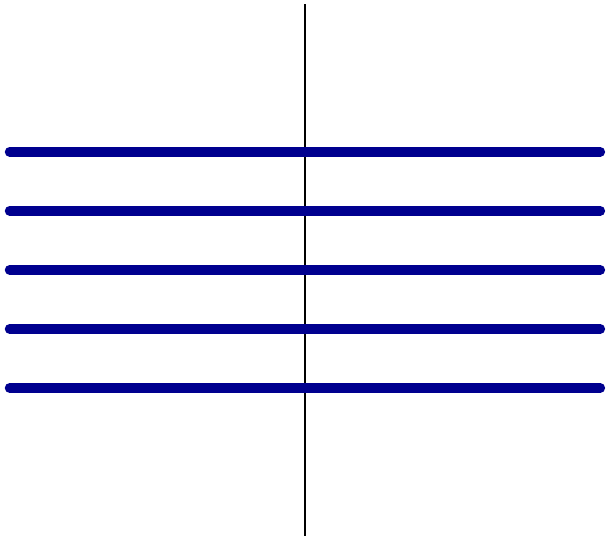}
\hspace*{5ex}
\includegraphics[width=16ex]{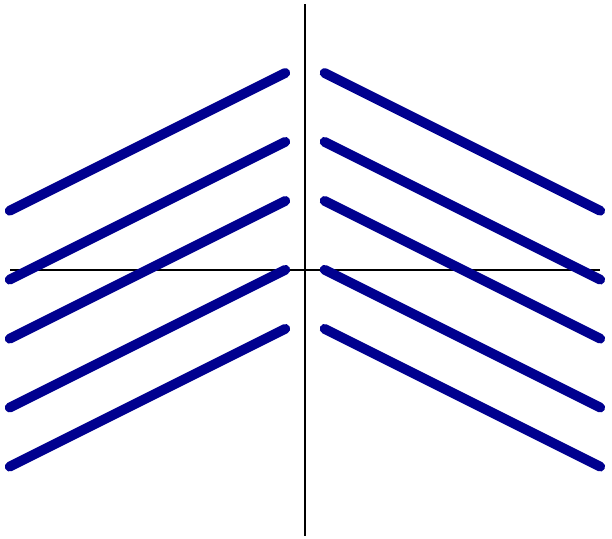}
\hspace*{5ex}
\includegraphics[width=16ex]{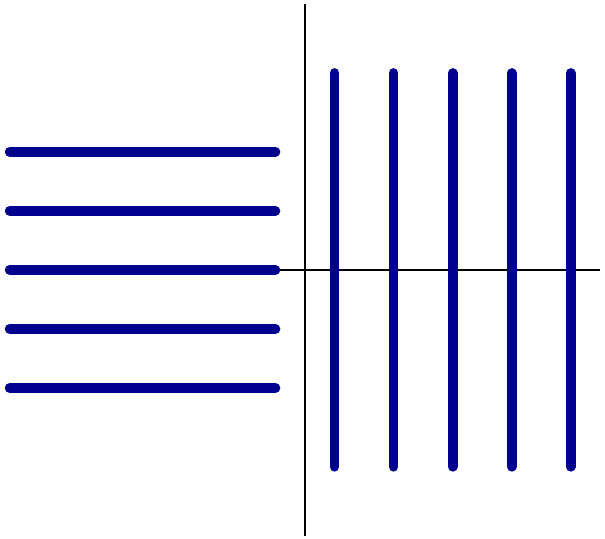}
\caption{From left to right, schematic plots of the projections on the horizontal plane of: convective rolls, a symmetric domain wall between two set of rolls rotated by opposite angles, and an orthogonal domain wall.}
  \label{f:plots}
\end{figure}

The B\'enard-Rayleigh convection problem is a classical problem in fluid mechanics. It concerns the flow of a three-dimensional viscous fluid layer situated between two horizontal parallel plates and heated from below. Upon increasing the difference of temperature between the two plates, the simple conduction state looses stability at a critical value of the temperature difference. In terms of nondimensional parameters this instability occurs at a critical value $\mathcal R_c$ of the Rayleigh number. Beyond the instability threshold, a convective regime develops in which patterns are formed, such as convective rolls, hexagons, or squares. Observed patterns are often accompanied by defects \cite{Bod,Man}. A particular class of such defects are domain walls which occur between rolls with different orientations (see Figure~\ref{f:plots}).

Mathematically, the governing equations are the Navier-Stokes equations coupled with an equation for the temperature, and completed by boundary conditions at the two plates (e.g., see \cite{Kosh}). Observed patterns are then found as particular steady solutions of these equations. Since the pioneering works of Yudovich \cite{Uk-Yu, Yu66,Yu67a,Yu67b}, Rabinowitz \cite{Rabinowitz}, and G\"{o}rtler et al \cite{GKS69} in the sixties, the existence of patterns was studied in various works by different authors (e.g., see \cite{GSS,Kosh,Br-Io} and the references therein).
Very recently, the existence of symmetric domain walls has been shown in \cite{HI21a,HI21b}.

Handling the full governing equations being often technically challenging, alternative studies rely on simpler amplitude equations which provide approximate descriptions of solutions in particular parameter regimes. We adopt this type of approach for the existence problem for orthogonal domain walls.

As a first step, we rigorously derive the system of amplitude equations \eqref{e:A0}-\eqref{e:B0} in the parameter regime of Rayleigh numbers slightly above the threshold of convective instability~$\mathcal R_c$. We apply the reduction procedure used in \cite{HI21a,HI21b} for the analysis of symmetric domain walls. Starting from a formulation of the steady governing equations as an infinite-dimensional dynamical system in which the horizontal coordinate $x$ plays the role of evolutionary variable, we apply a center manifold reduction and obtain a $12$-dimensional reduced dynamical system. Then, we compute a normal form for this reduced system and find the system \eqref{e:A0}-\eqref{e:B0} to leading order after a rescaling of the normal form. This first step is carried out in Section~\ref{s:system}.

Solutions of the system \eqref{e:A0}-\eqref{e:B0} provide leading order approximations of solutions of the full governing equations. In particular, the equilibrium $(A_0,B_0)=(0,1)$ of the system \eqref{e:A0}-\eqref{e:B0} gives an approximation of convection rolls bifurcating for Rayleigh numbers $\mathcal R>\mathcal R_c$ close to $\mathcal R_c$, whereas the equilibrium $(A_0,B_0)=(1,0)$ of the system \eqref{e:A0}-\eqref{e:B0} gives the same convection rolls but rotated by an angle $\pi/2$. A heteroclinic orbit connecting these two equilibria provides then an approximation of orthogonal domain walls. Our main result shows the existence of such heteroclinic orbits for the system \eqref{e:A0}-\eqref{e:B0}.
\begin{Theorem}\label{t:main}
  For any $\varepsilon>0$ and $g>1$, the system \eqref{e:A0}-\eqref{e:B0} possesses a smooth real-valued heteroclinic solution $(A_0,B_0)=(A_{\varepsilon,g},B_{\varepsilon,g})$ with the following properties:
  \begin{enumerate}
  \item $\lim_{x\to-\infty}(A_{\varepsilon,g}(x),B_{\varepsilon,g}(x))=(1,0)$ and
$\lim_{x\to\infty}(A_{\varepsilon,g}(x),B_{\varepsilon,g}(x))=(0,1)$;
    \item $B_{\varepsilon,g}(x)\geq 0$, for all $x\in\R$;
    \item for fixed $\varepsilon>0$, $\lim_{g\to 1^+} \sup_{x\in\R} |A_{\varepsilon,g}(x)^2+B_{\varepsilon,g}(x)^2-1|=0.$
      \end{enumerate}
\end{Theorem}

After some rescaling, the limit $\varepsilon=0$ is also considered, as it could give indications of how the heteroclinic orbits look like for small $\varepsilon>0$.

The proof of Theorem~\ref{t:main} is given in Section~\ref{s:heteroclinic}. The heteroclinic solution being real-valued, it is a solution of the $6$-dimensional system obtained by restricting the system \eqref{e:A0}-\eqref{e:B0} to real-valued functions $A_0$ and $B_0$. We use a variational method in which the heteroclinic solution is found as a minimizer of a rescaled functional. A compactness by concentration type argument is used to prove the convergence of minimizing sequences towards the heteroclinic solution. We point out that a complete existence proof for orthogonal domain walls for the Rayleigh-B\'enard convection problem requires an additional analysis consisting in showing that the heteroclinic orbit in Theorem~\ref{t:main} persists as a perturbed heteroclinic of the full $12$-dimensional reduced system, hence without restricting to the leading order system  \eqref{e:A0}-\eqref{e:B0}. Such a proof was given for symmetric domain walls in \cite{HI21a,HI21b}.

\medskip

\noindent
{\bf Acknowledgments:}
M.H. was partially supported by the project Optimal (ANR-20-CE30-0004) and the EUR EIPHI program (ANR-17-EURE-0002).

\section{Derivation of the amplitude equations}\label{s:system}

Relying upon a center manifold reduction and a normal forms analysis, we derive the system of amplitude equations \eqref{e:A0}-\eqref{e:B0} from the B\'enard-Rayleigh convection problem. This derivation being similar to that of the leading order systems in \cite{HI21a, HI21b}, we recall the main steps, focus on differences, and refer to these works for further details.

\subsection{Formulation of the hydrodynamic problem}

We consider the formulation as a dynamical system of the governing equations for the steady convection problem from \cite{HI21a}. In Cartesian coordinates $(x,y,z)\in \mathbb{R}^{3}$, where $(x,y)$ are the horizontal coordinates and $z$ is the vertical coordinate, after rescaling variables, the fluid occupies the domain $\mathbb{R}^{2}\times (0,1)$. The physical variables are the particle velocity $\mathbf{V}=(V_{x},V_{y},V_{z})$, the deviation of the temperature from the conduction profile $\theta $, and the pressure $p$. There are two dimensionless parameters, the Rayleigh number $\mathcal{R}$ and the Prandtl number $\mathcal{P}$.

Taking the horizontal coordinate $x$ as evolutionary variable, the governing equations are written as a system of the form
\begin{equation}
\partial _{x}\mathbf{U}=\mathcal{L}_{\mu }\mathbf{U}+\mathcal{B}_{\mu }( 
\mathbf{U},\mathbf{U}),  \label{newsyst}
\end{equation}
with $\mathbf{U}=(V_{x},V_{\bot},W_{x},W_{\bot},\theta ,\phi )$ an 8-components vector, in which $V_{\bot }=(V_{y},V_{z})$, $W_{\bot }=(W_{y},W_{z})$, and $\mathbf{W}=(W_{x},W_{\bot })$ and $\phi$ are additional variables defined by
\begin{equation}
\mathbf{W}=\mu ^{-1}\partial _{x}\mathbf{V}-p\mathbf{e}_{x},\quad \phi
=\partial _{x}\theta , \label{W}
\end{equation}
where $\mathbf{e}_{x}=(1,0,0)$.
The parameter $\mu$ is the square root of the Rayleigh number, $\mu =\mathcal{R}^{1/2}$, and $\mathcal{L}_{\mu }$ and $\mathcal{B}_{\mu }$ in the right hand side of \eqref{newsyst} are linear
and quadratic operators, respectively, defined by 
\begin{equation*}
\mathcal{L}_{\mu }\mathbf{U=}\left( 
\begin{array}{c}
-\nabla _{\bot }\cdot V_{\bot } \\[0.5ex] 
\mu W_{\bot } \\[0.5ex] 
-\mu ^{-1}\Delta _{\bot }V_{x} \\[0.5ex] 
-\mu ^{-1}\Delta _{\bot }V_{\bot }-\theta \mathbf{e}_{z}-\mu ^{-1}\nabla
_{\bot }(\nabla _{\bot }\cdot V_{\bot })-\nabla _{\bot }W_{x} \\[0.5ex] 
\phi \\[0.5ex] 
-\Delta _{\bot }\theta -\mu V_{z}%
\end{array}
\right) ,
\end{equation*}
\begin{equation*}
\mathcal{B}_{\mu }\mathbf{(U,U)}=\left( 
\begin{array}{c}
0 \\ 
0 \\ 
\mathcal{P}^{-1}\big((V_{\bot }\cdot \nabla _{\bot })V_{x}-V_{x}(\nabla
_{\bot }\cdot V_{\bot })\big) \\[0.5ex] 
\mathcal{P}^{-1}\big((V_{\bot }\cdot \nabla _{\bot })V_{\bot }+\mu
V_{x}W_{\bot }\big) \\[0.5ex] 
0 \\[0.5ex] 
\mu \big((V_{\bot }\cdot \nabla _{\bot })\theta +V_{x}\phi \big)%
\end{array}
\right) ,
\end{equation*}
where $\Delta_\bot=\partial_{yy}+\partial_{zz}$, $\nabla_{\bot }=(\partial_{y},\partial_{z})$, and here $\mathbf{e}_z=(0,1)$.

The phase space $\mathcal X$ for the dynamical system \eqref{newsyst} and the domain of definition $\mathcal Z$ of the linear operator $\mathcal L_\mu$ include the boundary conditions and a condition on the flux. We consider periodic boundary conditions in $y$ and the case of ``rigid-rigid'' boundary conditions in $z$:\footnote{The subsequent analysis remains valid for other types of boundary conditions in $z$; see \cite[Section 8]{HI21a} and \cite[Section 2]{HI21b} for the definition of the spaces $\mathcal X$ and $\mathcal Z$ in the cases of ``free-free'' and ``rigid-free'' boundary conditions, respectively.}
\begin{equation}
\mathbf{V}|_{z=0,1}=0,\quad \theta |_{z=0,1}=0.  \label{e:bcrigid}
\end{equation}
Taking the period $2\pi/k$ in $y$, for some fixed $k>0$, a direct calculation shows that the derivative of the flux 
\begin{equation*}
  \mathcal{F}(x)=\int_{\Omega _{per}}V_{x}\,dy\,dz,\quad
  \Omega _{per}=(0,2\pi/k_y)\times(0,1),
\end{equation*}
vanishes, hence $\mathcal F(x)$ is a constant function (see \cite[Section 3]{HI21a}). Fixing the constant flux to $0$, the phase space $\mathcal X$ is defined by
\begin{eqnarray*}
\mathcal{X}=\big\{\mathbf{U}\in (H_{per}^{1}(\Omega ))^{3}\times
(L_{per}^{2}(\Omega ))^{3}\times H_{per}^{1}(\Omega )\times
L_{per}^{2}(\Omega )\;;\; && \\
\hspace*{2ex}V_{x}=V_{\bot }=\theta =0\text{ on }z=0,1,\text{ and }%
\int_{\Omega _{per}}V_{x}\,dy\,dz=0\big\}, &&
\end{eqnarray*}
where the subscript $per$ means that the functions are $2\pi /k$-periodic in $y$
(for simplicity, we have written $V_x=V_{\bot}=\theta=0$ although these vectors do not have the same dimension).
The domain of definition $\mathcal{Z}$ of the linear operator $\mathcal{L}_{\mu }$ is defined by
\begin{eqnarray*}
\mathcal{Z}=\big\{\mathbf{U}\in \mathcal{X}\cap (H_{per}^{2}(\Omega
))^{3}\times (H_{per}^{1}(\Omega ))^{3}\times H_{per}^{2}(\Omega )\times
H_{per}^{1}(\Omega )\;;\; && \\[0.5ex]
\nabla _{\bot }\cdot V_{\bot }=W_{\bot }=\phi =0\text{ on }z=0,1\big\}, &&
\end{eqnarray*}
such that $\mathcal{L}_{\mu }$ is a closed operator. The space $\mathcal{Z}$ being compactly embedded in $\mathcal X$, the operator $\mathcal{L}_{\mu }$ has compact resolvent and therefore purely point spectrum consisting of isolated eigenvalues with finite algebraic multiplicity.

As a consequence of the symmetries of the hydrodynamic problem, the dynamical system \eqref{newsyst} is reversible with reversibility symmetry  
\begin{equation*}
\boldsymbol{S}_{1}\boldsymbol{U}(y,z) =(-V_{x},V_{\bot },W_{x},-W_{\bot },\theta
,-\phi )(y,z),\quad \boldsymbol{U}\in\mathcal{X},
\end{equation*}
which anti-commutes with $\mathcal{L}_{\mu }$ and $\mathcal{B}_{\mu }$, and $O(2)$-equivariant  with discrete symmetry 
\begin{equation*}
\boldsymbol{S}_{2}\boldsymbol{U}(y,z)=(V_{x},-V_{y},V_{z},W_{x},-W_{y},W_{z},
\theta ,\phi )(-y,z), \quad \boldsymbol{U}\in\mathcal{X},
\end{equation*}
and continuous symmetry $(\boldsymbol{\tau }_{a})_{a\in\mathbb{\ R}/2\pi\mathbb{Z}}$,
\begin{equation*}
  \boldsymbol{\tau}_{a}\boldsymbol{U}(y,z)=\boldsymbol{U}(y+a/k_y,z),
  \quad \boldsymbol{U}\in\mathcal{X},
\end{equation*}
which commute with $\mathcal{L}_{\mu }$ and $\mathcal{B}_{\mu }$ and satisfy
\begin{equation*}
\boldsymbol{\tau}_{a}\boldsymbol{S}_2=\boldsymbol{S}_2\boldsymbol{\tau}_{-a}, \quad 
\boldsymbol{\tau}_{0}=\boldsymbol{\tau} _{2\pi }=\mathbb{I}.
\end{equation*}
The symmetries $\boldsymbol S_1$ and $\boldsymbol S_2$ follow from the reflections $x\mapsto-x$ and $y\mapsto-y$, respectively, whereas the continuous symmetry $\boldsymbol{\tau}_{a}$ is a consequence of the invariance under translations in $y$ of the governing equations.\footnote{In the case of ``rigid-rigid'' and ``free-free'' boundary conditions, there is an additional vertical reflection symmetry $z\mapsto1-z$ leading to the symmetry 
$\boldsymbol{S}_{3}\mathbf{U}(y,z)=(V_{x},V_{y},-V_{z},W_{x},W_{y},-W_{z},-
\theta ,-\phi )(y,1-z)$
which commutes with $\mathcal{L}_{\mu }$ and $\mathcal{B}_{\mu }$. Aiming for a result which is also valid in the case of ``rigid-free'' boundary conditions, we do not make use of this symmetry.}
In addition, the system does not change when adding any constant to the new variable $W_{x}$, i.e., it is invariant under the action of the one-parameter family of maps $(\boldsymbol{T}_{b})_{b\in\mathbb{R}}$ defined by
\begin{equation}
\boldsymbol{T}_{b}\boldsymbol{U} = \boldsymbol{U}+b\boldsymbol{\varphi }_{0},\quad 
\boldsymbol{\varphi }_{0} =(0,0,0,1,0,0,0,0)^{t},
 \quad \boldsymbol{U}\in\mathcal{X}.  \label{phi0}
\end{equation}
 
In this setting, the classical convection rolls are equilibria of the dynamical system \eqref{newsyst}. As explained in \cite[Section 4]{HI21a}, these rolls provide a circle of equilibria $\boldsymbol{\tau }_{a}(\mathbf{U}_{k,\mu }^{* })$, for $a\in \mathbb{R}/2\pi \mathbb{Z}$, of the dynamical system \eqref{newsyst}
which bifurcate for $\mu >\mu _{0}(k)$ sufficiently close to a critical value $\mu _{0}(k)$, for any fixed wavenumber $k$. Due to the rotation invariance of the hydrodynamic problem, horizontally rotated rolls are solutions of the dynamical system \eqref{newsyst}. In particular, for the rotation angle of $\pi/2$, we obtain solutions which are $2\pi/k$-periodic in $x$ and constant in $y$. Orthogonal domain walls could then be constructed as heteroclinic orbits connecting these latter periodic solutions with the equilibria $\mathbf{U}_{k,\mu }^{* }$. According to the classical theory, the map $k\mapsto\mu_0(k)$ is analytical in $k$ and has a strict global minimum at $k=k_c$ where $\mu_0''(k_c)>0$. The values $k_c$, $\mu_0(k_c)$ and $\mu_0''(k_c)$ depend on the imposed boundary conditions at $z=0,1$ and can be computed numerically \cite[Section 2]{HI21a}.

\subsection{Reduced dynamics}

We consider the parameter regime with $(k,\mu)$ close to $(k_c,\mu_c)$, where $\mu_c=\mu_0(k_c)$. We set
\[
\mu=\mu_c+\tilde\mu,\quad k=k_c(1+\tilde k),
\]
in which $\tilde\mu$ and $\tilde k$ are small parameters.
We also eliminate the dependence on $k$ of the phase space $\mathcal X$ of the dynamical system \eqref{newsyst} by normalizing to $2\pi/k_c$ the period in $y$ of the solutions. The resulting system is of the form \eqref{newsyst} in which now $\Delta_\bot=(1+\tilde k)^2\partial_{yy}+\partial_{zz}$, $\nabla_{\bot }=((1+\tilde k)\partial_{y},\partial_{z})$, and its phase space is $\mathcal X$ with $k=k_c$. We write this system in the form
\begin{equation}
\partial _{x}\mathbf{U}=\mathcal{L}_{{c}}\mathbf{U}+\mathcal{R}(\mathbf{U},\tilde\mu,\tilde k),  \label{newsystb}
\end{equation}
where 
\begin{equation}\label{e:linear}
  \mathcal L_c = \mathcal L_{\mu_c}\big|_{\tilde k=0},\quad
\mathcal{R}(\mathbf{U},\tilde\mu,\tilde k)=(\mathcal{L}_{\mu }-\mathcal{L}_{\mu
_{c}}\big|_{\tilde k=0})\mathbf{U} +\mathcal{B}_{\mu }(\mathbf{U},\mathbf{U}),
\end{equation}
and $\mathcal R$ is a smooth map from $\mathcal{Z}\times (-\mu _{c},\infty )\times \R$ into $\mathcal{X }$ satisfying
\begin{equation}\label{e:nonlinear}
\mathcal{R}(0,\tilde\mu,\tilde k)=0,\quad D_{\mathbf{U}}\mathcal{R}(0,0,0)=0.
\end{equation}

We apply a center manifold theorem to obtain a reduced system of ordinary differential equations which describes the dynamics of \eqref{newsystb} in a neighborhood of the equilibrium $\mathbf U=0$ for small $(\tilde\mu,\tilde k)$.
The arguments are the same as the ones from \cite[Section 5]{HI21a}, except for the purely imaginary eigenvalues of the linear operator $\mathcal L_c$ which are different. The following result is obtained  by taking the limit $\alpha=0$ in the result from \cite[Lemma 4.2]{HI21a}.

\begin{Lemma}
  \label{l:spectrum}
The center spectrum of the linear operator $\mathcal{L}_{c}$ consists of the three eigenvalues $0,\pm ik_c$ with the following properties.

\begin{enumerate}
\item The eigenvalue $0$ has algebraic multiplicity $9$ and geometric multiplicity $3$, and the complex conjugated eigenvalues $\pm ik_c$ are algebraically double and geometrically simple.
\item For the eigenvalue $0$, there are three linearly independent eigenvectors: $\boldsymbol{\varphi }_{0}$ given by \eqref{phi0}, $\boldsymbol{\zeta}_0$ of the form $\boldsymbol{\zeta}_{0}(y,z)=\widehat{U}_{k_{c}}(z)e^{ik_{c}y}$, with $\widehat{U}_{k_{c}}(z)\in\C^8$, and the complex conjugated vector $\boldsymbol{\bar\zeta}_0$, and two chains of generalized eigenvectors: $\boldsymbol{\zeta}_{1}, \boldsymbol{\zeta}_{2}, \boldsymbol{\zeta}_{3}$ associated to $\boldsymbol{\zeta}_{0}$, \footnote{For our purposes, we do not need the explicit formulas for eigenvectors and generalized eigenvectors which can be obtained from \cite[Section 4]{HI21a}.} 
  \[
\mathcal L_c\boldsymbol{\zeta}_{1}=\boldsymbol{\zeta}_{0},\quad
\mathcal L_c\boldsymbol{\zeta}_{2}=\boldsymbol{\zeta}_{1},\quad
\mathcal L_c\boldsymbol{\zeta}_{3}=\boldsymbol{\zeta}_{2},
  \]
and the conjugated vectors  $\overline{\boldsymbol{\zeta}_{1}}, \overline{\boldsymbol{\zeta}_{2}}, \overline{\boldsymbol{\zeta}_{3}}$ associated to $\overline{\boldsymbol{\zeta}_{0}}$. The eigenvector $\boldsymbol{\varphi }_{0}$ is invariant under the actions of $\boldsymbol{S}_1$, $\boldsymbol{S}_{2}$, and $\boldsymbol{\tau}_{a}$, and the other generalized eigenvectors satisfy:
  \begin{eqnarray*}
&&\boldsymbol{S}_{1}\boldsymbol{\zeta }_{0}=\boldsymbol{\zeta }_{0},\quad \boldsymbol{S}_{2}\boldsymbol{\zeta }
_{0}=\overline{\boldsymbol{\zeta }_{0}},\quad \boldsymbol{\tau }_{a}\boldsymbol{\zeta }%
_{0}=e^{ia}\boldsymbol{\zeta }_{0}, \\
&&\boldsymbol{S}_{1}\boldsymbol{\zeta }_{1}=-\boldsymbol{\zeta }_{1},\quad \boldsymbol{S}_{2}\boldsymbol{\zeta 
}_{1}=\overline{\boldsymbol{\zeta }_{1}},\quad \boldsymbol{\tau }_{a}\boldsymbol{\zeta 
}_{1}=e^{ia}\boldsymbol{\zeta }_{1}, \\
&&\boldsymbol{S}_{1}\boldsymbol{\zeta }_{2}=\boldsymbol{\zeta }_{2},\quad \boldsymbol{S}_{2}\boldsymbol{\zeta }
_{2}=\overline{\boldsymbol{\zeta }_{2}},\quad \boldsymbol{\tau }_{a}\boldsymbol{\zeta }%
_{2}=e^{ia}\boldsymbol{\zeta }_{2}, \\
&&\boldsymbol{S}_{1}%
\boldsymbol{\zeta }_{3}=-\boldsymbol{\zeta }_{3},\quad \boldsymbol{S}_{2}\boldsymbol{\zeta 
}_{3}=\overline{\boldsymbol{\zeta }_{3}},\quad \boldsymbol{\tau }_{a}\boldsymbol{\zeta 
}_{3}=e^{ia}\boldsymbol{\zeta }_{3}.
\end{eqnarray*}

\item For the eigenvalue $ik_c$, there is one eigenvector $\boldsymbol{\xi }_{0}$ of the form $\boldsymbol{\xi }_{0}(y,z)=\widehat{\mathbf{U}}_{0}(z)\in\C^8$,
and an associated generalized eigenvector $\boldsymbol{\xi }_{1}$ with the properties
\[
(\mathcal{L}_{{c}}-ik_{c})\boldsymbol{\xi}_{1}=\boldsymbol{\xi}_{0},
\]
and
\begin{eqnarray*}
&&\boldsymbol{S}_{1}\boldsymbol{\xi }_{0}=\overline{\boldsymbol{\xi }_{0}}
,\quad \boldsymbol{S}_{2}\boldsymbol{\xi }_{0}=\boldsymbol{\xi }_{0},\quad 
\boldsymbol{\tau }_{a}\boldsymbol{\xi }_{0}=\boldsymbol{\xi }_{0}, \\
&&\boldsymbol{S}_{1}\boldsymbol\xi_1=-\overline{\boldsymbol\xi_1}
,\quad \boldsymbol{S}_{2}\boldsymbol\xi_1=\boldsymbol\xi_1,\quad \boldsymbol{\tau }_{a}%
\boldsymbol{\xi }_{1}=\boldsymbol\xi_1.
\end{eqnarray*}
The complex conjugated vectors $\overline{\boldsymbol{\xi }_{0}}$ and $\overline{\boldsymbol{\xi }_{1}}$ are eigenvector and generalized eigenvector, respectively, for the eigenvalue $-ik_{c}$.
\end{enumerate}
\end{Lemma}

As a result of the center manifold theorem, we obtain that the small bounded solutions of the infinite-dimensional dynamical system \eqref{newsystb} belong to a $13$-dimensional center manifold, for any sufficiently small $\tilde \mu$ and $\tilde k$, and are of the form
\begin{eqnarray*} \label{centermanifold}
  \boldsymbol{U} &=&w\boldsymbol{\phi }_{0}
  +A_{0}\boldsymbol{\zeta }_{0}+A_{1}\boldsymbol{\zeta }_{1}+A_{2}\boldsymbol{\zeta }_{2}+A_{3}\boldsymbol{\zeta }_{3}+B_{0}\boldsymbol{\xi }_{0}+B_{1}\boldsymbol{\xi }_{1}\\
  && +\overline{A_{0}\boldsymbol{\zeta }_{0}}+\overline{A_{1}\boldsymbol{\zeta }_{1}}+\overline{A_{2}\boldsymbol{\zeta }_{2}}+\overline{A_{3}\boldsymbol{\zeta }_{3}}+\overline{B_{0}\boldsymbol{\xi }_{0}}+\overline{B_{1}\boldsymbol{\xi }_{1}}
  +\boldsymbol{\Phi }(X,\overline{X},\widetilde{\mu },\widetilde{k}), 
\end{eqnarray*}
in which $w$ and $X=(A_{0},A_{1},A_{2},A_{3},B_{0},B_{1})$ are $x$-dependent functions with values in $\R$ and $\C^6$, respectively, and $\Phi$ is of class $C^m$ in its arguments, for any fixed $m\geq1$. The eigenvectors and generalized eigenvalues being complex-valued, it is convenient to use here complex variables $(X,\overline X)$, instead of $12$ real variables, hence by identifying $\R^{12}$ with the space $\C^6\times\overline{\C^6}=\{(Z,\overline Z)\;;\; Z\in\C^6\}$.

The reduced $13$-dimensional system for $w$, $X$, and $\overline X$ inherits the properties of the infinite-dimensional dynamical system \eqref{newsyst}. In particular, the invariance of \eqref{newsyst} under the action of $\boldsymbol T_b$, implies that the reduced vector field is invariant under the action of the induced transformation $w\mapsto w+b$, for any $b\in\R$, and therefore does not depend on $w$. Consequently, the equations for $w$ and $(X,\overline X)$ in the reduced system are decoupled,
\begin{eqnarray*}
  \frac{dw}{dx} =h(X,\overline{X},\widetilde{\mu },\widetilde{k}),  \end{eqnarray*}
and
\begin{eqnarray}
  \frac{dX}{dx} =F(X,\overline{X},\widetilde{\mu },\widetilde{k}),
\quad 
\frac{d\overline X}{dx} =\overline{F(X,\overline{X},\widetilde{\mu },\widetilde{k})}.
  \label{reduced} 
\end{eqnarray}
Taking into account the properties \eqref{e:linear}-\eqref{e:nonlinear} and the result in Lemma~\ref{l:spectrum} we obtain that
\begin{equation}\label{e:propred}
F(0,0,\tilde\mu,\tilde k)=0, \quad
D_{X}F(0,0,0,0)=L_c,\quad D_{\overline X}F(0,0,0,0)=0,
\end{equation}
where $L_c$ is the $6\times6$ Jordan matrix  
\begin{equation} \label{e:L0}
L_c=\left( 
\begin{array}{cc}
L_0 & 0 \\ 
0 & L_1 \\
\end{array}
\right),\quad
L_0=\left( 
\begin{array}{cccc}
0 & 1 & 0 & 0 \\ 
0 & 0 & 1 & 0 \\
0 & 0 & 0 & 1 \\ 
0 & 0 & 0 & 0
\end{array}
\right),\quad
L_1=\left( 
\begin{array}{cc}
ik_c & 1 \\[0.5ex] 
0 & ik_c 
\end{array}
\right).
\end{equation}
In addition, from the symmetry properties of the eigenvectors and generalized eigenvectors in Lemma~\ref{l:spectrum}, we deduce their actions on the variable $X$,
\begin{eqnarray}\label{e:S1}
&&\boldsymbol{S}_{1}(A_{0},A_{1},A_{2},A_{3},B_{0},B_{1})
=(A_{0},-A_{1},A_{2},-A_{3},\overline{B_{0}},-\overline{B_{1}}),
 \\ \label{e:S2}
&&\boldsymbol{S}_{2}(A_{0},A_{1},A_{2},A_{3},B_{0},B_{1})
 =(\overline{A_{0}},
 \overline{A_{1}},\overline{A_{2}},\overline{A_{3}},B_{0},B_{1}),
 \\ \label{e:tau}
&&\boldsymbol{\tau }_{a}(A_{0},A_{1},A_{2},A_{3},B_{0},B_{1})
=(e^{ia}A_{0},e^{ia}A_{1},e^{ia}A_{2},e^{ia}A_{3},B_{0},B_{1}).
\end{eqnarray}
Then, the vector field in the reduced system \eqref{reduced} anti-commutes with $\boldsymbol S_1$ and commutes with $\boldsymbol S_2$ and $\boldsymbol\tau_a$.
Notice that the equivariance under the action of $\boldsymbol S_2$ implies that the reduced system leaves invariant the $8$-dimensional subspace $\{(X,\overline{X})\;;\;A_j=\overline{A_j},\; j=0,1,2,3\}$. Solutions in this subspace correspond to solutions of \eqref{newsyst} which are even in $y$.  
There is a second invariant subspace $\{(X,\overline{X})\;;\;A_j=0,\; j=0,1,2,3\}$, which corresponds to solutions of \eqref{newsyst} which do not depend on $y$.

\subsection{Leading order dynamics}

Next, we obtain a cubic normal form for the reduced system \eqref{reduced}. The following result holds for general $12$-dimensional vector fields which satisfy \eqref{e:propred} and have the symmetries \eqref{e:S1}-\eqref{e:tau}.

\begin{Lemma}\label{l:normalform}
Consider a system of ordinary differential equations of the form \eqref{reduced} in which the vector field $F$ is of class $C^m$, for some $m\geq4$, in a neighborhood $\mathcal{U}_{1}\times \overline{\mathcal{U}_{1}}\times \mathcal{U}_{2}\subset\mathbb{C}^{6}\times\overline{\mathbb{C}^6}\times\R^2$ of the origin. Assume that the properties \eqref{e:propred} hold and that $F$ anti-commutes with $\boldsymbol S_1$ in \eqref{e:S1} and commutes with $\boldsymbol S_2$ in \eqref{e:S2} and $\boldsymbol\tau_a$ in \eqref{e:tau}.

There exist neighborhoods 
$\mathcal{V}_{1}$ and $\mathcal{V}_{2}$ of $0$ in $\mathbb{C}^{6}$ and $\mathbb{R}^2$, respectively, such that for any $(\tilde\mu,\tilde k)\ \in \mathcal{V}_{2} $, there is a polynomial ${P}(\cdot,\cdot,\tilde\mu,\tilde k): \mathbb{C}^{6}\times \overline{\mathbb{C}^6}\to\mathbb{C}^{6}$ of degree $3$ in the variables $(Z,\overline Z)$, such that for $Z\in \mathcal{V}_{1}$, the
change of variable 
\begin{equation*}  
X=Z+P(Z,\overline Z,\tilde\mu,\tilde k),
\end{equation*}
transforms the equation \eqref{reduced} into the normal form 
\begin{equation}
\frac{dZ}{dx}={L}_cZ+{N}(Z,\overline Z,\tilde\mu,\tilde k)+\rho (Z,\overline
Z,\tilde\mu,\tilde k),  \label{normal form}
\end{equation}
with the following properties:

\begin{enumerate}
\item the map $\rho $ belongs to $\mathcal{C}^{m}(\mathcal{V}_{1}\times 
\overline{\mathcal{V}_{1}}\times \mathcal{V}_{2}, \mathbb{C}^{6})$, and 
\begin{equation*}
\rho (Z,\overline Z,\tilde\mu,\tilde k
)=O(|(\tilde\mu,\tilde k)|^2\|Z\|+|(\tilde\mu,\tilde k)|\,\|Z\|^{3}+\|Z\|^{4});
\end{equation*}

\item both $N(\cdot,\cdot,\tilde\mu,\tilde k)$ and $\rho(\cdot,\cdot,\tilde\mu,\tilde k)$ anti-commute with $\boldsymbol{S}_1$ and commute with $\boldsymbol{S}_2$ and $\boldsymbol{\tau}_a$, for any $(\tilde\mu,\tilde k)\in \mathcal{V}_2$;

\item the six components $(N_0,N_1,N_2,N_3,M_0,M_1)$ of $N$ are of the form
  \begin{eqnarray*}
&&N_0 =iA_{0}P_{0},
  \\
&&N_1 =iA_1P_{0}+A_0P_1+b_7u_7,
\\
&&N_2=iA_{2}P_{0}+A_1P_1+iA_0P_2+ b_7v_7+c_8u_8+c_{9}u_{9} ,
\\
&&N_3=iA_{3}P_{0}+A_2P_1+iA_1P_2+A_0P_3+ b_7w_7+c_8v_8+c_{9}v_{9} +d_7u_7+d_{10}u_{10}+d_{11}u_{11} ,
\\
&&M_0=iB_0Q_0+\alpha_{12}u_{12},
\\
&&M_1=iB_1Q_0+B_0Q_1+\alpha_{12}v_{12}+i\beta_{12}u_{12}+i\beta_{13}u_{13}, 
\end{eqnarray*}
with
\begin{eqnarray*}
&&P_{0} =a_2u_2+a_4u_4, \\
&&P_1 =b_0 \tilde\mu+b_0'\tilde k +b_1u_1+b_3u_3+b_5u_5+b_6u_6, \\
&&P_2 =c_2u_2+c_4u_4, \\
&&P_3 =d_0 \tilde\mu+d_0'\tilde k +d_1u_1+d_3u_3+d_5u_5+d_6u_6, \\
&&Q_0 =\alpha_0 \tilde\mu+\alpha_0'\tilde k +\alpha_1u_1+\alpha_3u_3+\alpha_5u_5+\alpha_6u_6, \\
&&Q_1 =\beta_0 \tilde\mu+\beta_0'\tilde k +\beta_1u_1+\beta_3u_3+\beta_5u_5+\beta_6u_6,
\end{eqnarray*}
where $(A_0,A_1,A_2,A_3,A_4,B_0,B_1)$ are the six components of $Z$, the
coefficients $a_j$, $b_j$, $c_j$, $d_j$, $\alpha_j$, and $\beta _{j}$ are all real, and
\begin{eqnarray*}
&u_{1} =A_{0}\overline{A_{0}},\quad
 u_{2}=i(A_{0}\overline{A_1}-\overline{A_{0}}A_1),\quad
 u_{3}=A_{0}\overline{A_2}+\overline{A_{0}}A_2 -A_{1}\overline{A_{1}}, &\\
&u_{4}=i(A_{0}\overline{A_3}-\overline{A_{0}}A_3 -A_{1}\overline{A_2}+\overline{A_{1}}A_2),\quad
 u_{5} =B_0\overline{B_0},\quad
 u_{6}=i(B_{0}\overline{B_1}-\overline{B_{0}}B_1),&
   \end{eqnarray*}
  \begin{eqnarray*}
 &u_7= \overline{A_{0}} (A_1^2-2A_0A_2),\quad
 v_7= \overline{A_{1}} (A_1^2-2A_0A_2),\quad
 w_7= \overline{A_{2}} (A_1^2-2A_0A_2), &\\
 &
\! u_8= A_0v_3-A_1u_3, \quad
 v_8= A_1v_3-2A_2u_3, \quad
 v_3=\frac12 (3A_{0}\overline{A_3}+3\overline{A_{0}}A_3
 -A_{1}\overline{A_2}-\overline{A_{1}}A_2), \\
 &u_{9}= \frac12 A_0(B_{0}\overline{B_1}+\overline{B_{0}}B_1)-A_1u_5 ,\quad
 v_{9}=\frac12 A_0B_{1}\overline{B_1} -A_2u_5, &\\
& u_{10}= 3iA_3u_2
 +2A_2(A_{0}\overline{A_2}-\overline{A_{0}}A_2)
 -A_1(A_{1}\overline{A_2}-\overline{A_{1}}A_2) ,&\\
 & u_{11}= \frac12 A_0B_{1}\overline{B_1} +A_2u_5
 - \frac12 A_1(B_{0}\overline{B_1}+\overline{B_{0}}B_1) ,&
   \end{eqnarray*}
  \begin{eqnarray*}
 &u_{12}= \frac12 B_0(A_{0}\overline{A_1}+\overline{A_{0}}A_1) -
 B_1u_1 ,\quad v_{12}= B_0\overline{A_{1}}A_1- \frac12 B_1(A_{0}\overline{A_1}+\overline{A_{0}}A_1),&\\
 & u_{13}=B_0v_3-B_1u_3 .&  
\end{eqnarray*}
\end{enumerate}
\end{Lemma}

The proof of this result closely follows the arguments in the proofs from \cite[Lemma 6.1]{HI21a} and \cite[Theorem 2]{HI21b}. Differences are at computational level, only, and we therefore skip this proof. Also, similarly to \cite{HI21a} and \cite{HI21b} we can determine the coefficients of the leading order terms which will appear in the amplitude equations. We obtain that
\[
d_0'=\alpha_0'=\beta_0'=0,
\quad
d_0=-4k_c^2\beta_0>0,\quad d_1=-4k_c^2\beta_5<0,\quad
\frac{\beta _{1}}{\beta _{5}}=\frac{d_{5}}{d_{1}}:=g>0.  
\]
Following \cite[Section 6.3]{HI21a} and \cite[Section 4]{HI21b}, we assume that $\widetilde\mu>0$ and rescale variables
\begin{eqnarray*}
&x =\frac{1}{2\varepsilon k_{c}}\widetilde{x}, \quad
\widetilde{\mu } =\frac{4k_{c}^{2}}{-\beta _{0}}\varepsilon ^{4},\quad
\widetilde{k}=\varepsilon ^{2}\widehat{k},
&\\
&{A_{0}(x) =\frac{2k_{c}}{\sqrt{\beta _{5}}}\varepsilon ^{2}\widetilde{A_{0}}(\widetilde x),
\quad A_{1}(x)=\frac{4k_{c}^{2}}{\sqrt{\beta _{5}}}\varepsilon ^{3}\widetilde{A_{1}}(\widetilde x), }&\\
&\quad A_{2}(x)=\frac{8k_{c}^{3}}{\sqrt{\beta _{5}}}\varepsilon ^{4}\widetilde{A_{2}}(\widetilde x),
\quad A_{3}(x)=\frac{16k_{c}^{4}}{\sqrt{\beta _{5}}}\varepsilon ^{5}\widetilde{A_{3}}(\widetilde x), &\\
&B_{0}(x) =\frac{2k_{c}}{\sqrt{\beta _{5}}}\varepsilon ^{2}
e^{\frac i{2\varepsilon}\tilde x}\widetilde{B_{0}}(\widetilde x),
\quad B_{1}(x)=\frac{4k_{c}^{2}}{\sqrt{\beta _{5}}}\varepsilon ^{3}
e^{\frac i{2\varepsilon}\tilde x}\widetilde{B_{1}}(\widetilde x).&
\end{eqnarray*}
Notice here the exponential factor $e^{\frac i{2\varepsilon}\tilde x}$ in the formulas for $B_0$ and $B_1$.
Then, taking into account the properties of the coefficients above we obtain the rescaled system 
\begin{eqnarray*}
&&\frac{d\widetilde{A_{0}}}{d\widetilde{x}} =\widetilde{A_{1}}+O(|\varepsilon|^2(|\widehat k|^2+|\varepsilon|^2)),   \\
&&\frac{d\widetilde{A_{1}}}{d\widetilde{x}} =\widetilde{A_{2}}+O(|\widehat k|+|\varepsilon|^2),\\
&&\frac{d\widetilde{A_{2}}}{d\widetilde{x}} =\widetilde{A_{3}}+O(|\widehat k|+|\varepsilon|^2),   \\
&&\frac{d\widetilde{A_{3}}}{d\widetilde{x}} =
  \widetilde{A_{0}}(1-|\widetilde{A_{0}}|^{2}-g|\widetilde{B_{0}}|^{2})
    +O(|\widehat k|+|\varepsilon|),
\\
&&\frac{d\widetilde{B_{0}}}{d\widetilde{x}} = \widetilde{B_{1}}+O(|\varepsilon|(|\widehat k|+|\varepsilon|^2)),   \\
&&\frac{d\widetilde{B_{1}}}{d\widetilde{x}} =
\varepsilon ^{2}\widetilde{B_{0}}(-1+g|\widetilde{A_{0}}|^{2}+|\widetilde{B_{0}}|^{2}) +O(|\varepsilon|^2(|\widehat k|+|\varepsilon|)).  
\end{eqnarray*}
Keeping only the leading order terms in each equation, the resulting system is equivalent to the system \eqref{e:A0}-\eqref{e:B0}.

We point out that the real equilibrium $M_+=(0,1)$ of the system \eqref{e:A0}-\eqref{e:B0} corresponds to the roll solution $\mathbf U^*_{k,\mu}$ of the dynamical system \eqref{newsyst}, whereas the real equilibrium $M_-(1,0)$ corresponds to the same roll solution rotated by an angle $\pi/2$. Consequently, a domain wall connecting these two orthogonal rolls in the B\'enard-Rayleigh problem corresponds to a heteroclinic connection between these two real equilibria of the system  \eqref{e:A0}-\eqref{e:B0} (for further details, see \cite[Section 6.3]{HI21a} and \cite[Section 4.2]{HI21b}).

\section{Existence of a heteroclinic orbit}\label{s:heteroclinic}

In this section we prove the result in Theorem~\ref{t:main}. We restrict to real-valued solutions and choose new scales by taking
\begin{equation*}
\epsilon=\varepsilon ^{4}>0,\quad
\bar{x}=\epsilon ^{1/4}x,\quad
A_{0}(x)=\bar{A}(\bar{x}),\quad
B_{0}(x)=\bar{B}(\bar{x}).
\end{equation*}
Then, the system \eqref{e:A0}-\eqref{e:B0} becomes, after suppression of bars,
\begin{eqnarray}
&&\epsilon \frac{d^{4}A}{dx^{4}} =A(1-A^{2}-gB^{2}),
\label{e:A} \\
&&\frac{d^{2}B}{dx^{2}} =B(-1+gA^{2}+B^{2}).
\label{e:B}
\end{eqnarray}
We construct the heteroclinic orbit as a minimizer of the functional
\begin{equation*}
J_{\epsilon }(A,B):=\int_{\mathbb{R}}\Big(\frac{\epsilon }{2}A^{\prime
\prime 2}+\frac{1}{2}B^{\prime }{}^{2}+\frac{1}{4}(A^{2}+B^{2}-1)^{2}+\frac{1
}{2}(g-1)A^{2}B^{2}\Big)dx ,
\end{equation*}
on the set $X$ of real-valued functions $(A,B)\in H_{loc}^{2}(\mathbb{R})\times H_{loc}^{1}(\mathbb{R})$ such that 
\begin{equation} \label{e:limits}
\lim_{x\rightarrow -\infty }(A(x),B(x))=(1,0)~~\text{ and }
~~\lim_{x\rightarrow \infty }(A(x),B(x))=(0,1).
\end{equation}
For any $\epsilon>0$ and $g>1$ this functional is nonnegative, $J_{\epsilon }(A,B)\in[0,\infty]$. In fact $J_\epsilon$ is more generally defined on $H_{loc}^{2}(\mathbb{R})\times H_{loc}^{1}(\mathbb{R})$ with values in $[0,+\infty]$. 
  A delicate issue will be, once a solution $(A,B)\in H_{loc}^{2}(\mathbb{R})\times H_{loc}^{1}(\mathbb{R})$ is obtained with  $J_\epsilon(A,B)<\infty$, to check that indeed \eqref{e:limits} is satisfied.

Setting 
\begin{equation*}
P(A,B)=\frac{1}{4}\Big(A^{2}+B^{2}-1\Big)^{2}+\frac{1}{2}(g-1)A^{2}B^{2}\,,
\end{equation*}
a stationary point $(A,B)\in X$ of $J_{\epsilon }$ satisfies the system 
\begin{equation*}
\epsilon A^{\prime \prime \prime \prime }+\partial_AP(A,B)=0,\quad -B^{\prime \prime
}+\partial_BP(A,B)=0,
\end{equation*}
in the sense of distributions. This system is precisely the system \eqref{e:A}-\eqref{e:B}. Notice that a standard bootstrap argument shows that any solution 
$(A,B)\in H_{loc}^{2}(\mathbb{R})\times H_{loc}^{1}(\mathbb{R})$ is smooth if~$\epsilon>0$.

\subsection{The case $\boldsymbol{\epsilon=0}$}

Although the case $\epsilon=0$ is not part of the statement of Theorem~\ref{t:main}, we nevertheless mention this case, because it could give an additional insight on the problem when $\epsilon>0$ is small.

For $\epsilon= 0$ and $g>1$, we have the functional 
\begin{equation*}
J_{0}(A,B)=\int_{\mathbb{R}}\Big(\frac{1}{2}B^{\prime }{}^{2}+\frac{1}{4}
(A^{2}+B^{2}-1)^{2}+\frac{1}{2}(g-1)A^{2}B^{2}\Big)dx \in \lbrack 0,\infty ],
\end{equation*}
For fixed $B$, one can minimize with respect to $A$. Differentiating the map 
\begin{equation*}
A\rightarrow f(A):=\frac{1}{4}(A^{2}+B^{2}-1)^{2}+\frac{1}{2}
(g-1)A^{2}B^{2}\,,
\end{equation*}
one gets the equation for $A$: 
\begin{equation*}
(A^{2}+gB^{2}-1)A=0.
\end{equation*}
Hence critical points satisfy $A=0$ or $A^{2}=1-gB^{2}$ if $1-gB^{2}\geq0$. As $f^{\prime \prime}(A)=3A^{2}-(1-gB^{2})$, we see that if $1-gB^{2}>0$, then $f^{\prime \prime}(0)<0$ and the minimum of $f$ is reached at $A=\pm \sqrt{1-gB^{2}}$. Consequently, $A=0$ if $1-gB^{2}\leq 0$, and $A=\pm \sqrt{1-gB^{2}}$ if $1-gB^{2}\geq 0$, or equivalently,
\[
A^{2}=\max \{0,1-gB^{2}\}=(1-gB^{2})_{+}.
\]
Substituting $A^2$ above in $J_{0}(A,B)$, one gets the reduced functional 
\begin{equation*}
  J_{\text{red}}(B)=\int_{\mathbb{R}}\Big(\frac{1}{2}B^{\prime }{}^{2}
  +\frac{1}{4}((1-gB^{2})_{+}+B^{2}-1)^{2}+\frac{1}{2}(g-1)(1-gB^{2})_{+}B^{2}\Big)
dx\in \lbrack 0,\infty ],
\end{equation*}
which depends on $B$, only.

A stationary point $(A,B)\in C(\mathbb{R})\times H_{loc}^{1}(\mathbb{R})$ of $J_{0}$ satisfies 
\begin{equation*}
  \partial_{A}P(A,B)=0,~~-B^{\prime \prime }+\partial_{B}P(A,B)=0,
\end{equation*}
or equivalently,
\begin{equation*}
A(A^{2}+gB^{2}-1)=0,\quad -B^{\prime \prime }+B(gA^{2}+B^{2}-1)=0,
\end{equation*}
in the sense of distributions for the second equation, together with the property \eqref{e:limits} for the limits at $x=\pm\infty$.
Consequently, $A=0$ or $A^{2}=1-gB^{2}$ if $B^{2}\leq 1/g$, hence leading to  the equation for $B$ 
\begin{equation}\label{e:eqB}
B^{\prime \prime }=\left\{ 
\begin{array}{ll}
-B+B^{3} & \text{ if }B^{2}\geq 1/g  , \\ 
(g-1)B+(1-g^{2})B^{3} & \text{ if }B^{2}\leq 1/g .
\end{array}
\right.
\end{equation}
Observe that the right-hand side is continuous.
This problem has an increasing solution $B>0$ of class $C^2$ such that 
\begin{equation*}
\lim_{x\rightarrow -\infty }B(x)=0~\text{ and }~~\lim_{x\rightarrow \infty
}B(x)=1,
\end{equation*}
which gives a solution of \eqref{e:A}-\eqref{e:B} by taking  $A=\sqrt{1-gB^{2}}$ if $0<B \leq 1/\sqrt g$ and $A=0$ if $B\geq 1/\sqrt g$.

Indeed $B=1$ is an hyperbolic equilibrium of the first equation in \eqref{e:eqB} and $B=0$ is an hyperbolic equilibrium of the second equation in \eqref{e:eqB} because $g>1$. On the other hand, the first equation
possesses the invariant $|B^{\prime }|^{2}+B^{2}-\frac{1}{2}B^{4}$, which is 
$1/2$ at the equilibrium $B=1$, and the second equation possesses the invariant $|B^{\prime}|^{2}+(1-g)B^{2}-(1-g^{2})\frac{1}{2}B^{4}$, which is $0$ at the
equilibrium $B=0$. Let us study the curves for $B=\pm \sqrt{1/g}$ in the
plane $(B,B^{\prime })$. From $|B^{\prime }|^{2}+\frac{1}{g}-\frac{1}{2g^{2}}=1/2$ that corresponds to the first equation, one gets $|B^{\prime }|^{2}=\frac{1}{2}-\frac{1}{g}+\frac{1}{2g^{2}}$. This is the same value of $|B^{\prime }|^{2}$ that one gets by solving the second equation: $|B^{\prime }|^{2}+\frac{1-g}{g}-\frac{1-g^{2}}{2g^{2}}=0$.
This shows that $B^{\prime }$ is continuous (if its sign does not jump) at the junction of the two curves in the $(B,B^{\prime })$ plane. Hence there is an heteroclinic solution coming from $(B,B')=(0,0)$, staying on the set $|B^{\prime}|^{2}+(1-g)B^{2}-(1-g^{2})\frac{1}{2}B^{4}=0$ for $B\in[0,1/\sqrt g]$, then on the set $|B^{\prime }|^{2}+B^{2}-\frac{1}{2}B^{4}=1/2$ for $B\in[1/\sqrt g,1]$, and finally tending to $(B,B')=(1,0)$, thus providing a solution $(A,B)\in C(\R)\times C^2(\R)$ of our problem.

\subsection{Estimates}

From now on we assume that $\epsilon>0$ and $g>1$. 
Let $I$ be a closed interval of length $1$, $M_{1},M_{2}\in [0,\infty)$, and $(A,B)\in H^{2}(\overset{\circ }{I})\times H^{1}(\overset{\circ }{I})$ such that 
\begin{equation*}
\int_{I}\frac{1}{2}|A^{\prime \prime }(x)|^{2}dx\leq M_{1}~~\text{ and }
~~\int_{I}P(A(x),B(x))dx\leq M_{2}\,.
\end{equation*}
In particular, $A$, $B$, and the derivative $A^{\prime }$ are continuous functions. For all $x_{1}<x_{2}$ in $I$, one has 
\begin{equation*}
A(x_{2})-A(x_{1})=A^{\prime
}(x_{1})(x_{2}-x_{1})+\int_{x_{1}}^{x_{2}}(x_{2}-s)A^{\prime \prime }(s)ds
\end{equation*}
and thus 
\begin{eqnarray*}
  |A(x_{2})-A(x_{1})-A^{\prime }(x_{1})(x_{2}-x_{1})|&\leq&
  \left(\int_{x_{1}}^{x_{2}}(x_{2}-s)^{2}ds\right) ^{1/2} \cdot
  \left(\int_{x_{1}}^{x_{2}}A^{\prime \prime }(s)^{2}ds\right) ^{1/2}\\
&\leq&
3^{-1/2}|x_{2}-x_{1}|^{3/2}\sqrt{2M_{1}}.
\end{eqnarray*}
This remains true if $x_{2}\leq x_{1}$ in $I$.

As $\int_IP(A(x),B(x))dx\leq M_2$ and the integrand is nonnegative, there
exists $x_1 \in I$ such that $P(A(x_1),B(x_1))\leq M_2$. Observe that 
\begin{equation*}
P(A,B)\geq
K\min\{(B\pm 1)^2+A^2,\ B^2+(A\pm 1)^2\},
\end{equation*}
for some constant $K>0$, the Hessian 
\begin{equation*}
P^{\prime\prime}(A,B)=
\begin{pmatrix}
3A^2+gB^2-1 & 2gAB \\ 
2gAB & gA^2+3B^2-1
\end{pmatrix}
\end{equation*}
being positive definite at $(A,B)\in\{(\pm 1,0),(0,\pm 1)\}$ and the growth
of $P$ being quartic at infinity. Therefore 
\begin{equation*}
(\|(A,B)\|-1)^2\leq
\min\{(B\pm 1)^2+A^2,\ B^2+(A\pm 1)^2\} \leq 
\frac{P(A,B)} K,
\end{equation*}
\begin{equation*}
\|(A,B)\|\leq 1+ \sqrt{P(A,B)/K} ~~\text{ and thus }~~
\|(A(x_1),B(x_1))\|\leq 1+\sqrt {M_2/K}.
\end{equation*}
Let us check that  
$|A^{\prime}(x_1)|\leq 8\left(\sqrt{2M_1/3} +2+2\sqrt{M_2/K}\right)$.

After a possible translation in $x$, let us assume that $I=[-1/2,1/2]$. Let us also assume to be in the case $x_1\in[-1/2,0]$ and $A^{\prime}(x_1)\geq 0$. For $x\geq 1/4$, one has 
\begin{equation*}
A(x)\geq A(x_1)+A^{\prime}(x_1)(1/4)-\sqrt{2M_1/3}.
\end{equation*}
If moreover $A^{\prime}(x_1)(1/4)\geq 2\left( \sqrt{2M_1/3}+2+\sqrt{M_2/K}
\right)$, then $A(x)\geq 1$ for $x\in[1/4,1/2]$ and 
\begin{equation*}
\int_I P(A,B)dx \geq \int_{1/4}^{1/2}K(A-1)^2dx \geq
\int_{1/4}^{1/2}K\left(-|A(x_1)|+A^{\prime}(x_1)(1/4)-\sqrt{2M_1/3}
-1\right)^2dx
\end{equation*}
\begin{equation*}
\geq \int_{1/4}^{1/2}KA^{\prime}(x_1)^2 (1/{64})dx \geq
KA^{\prime}(x_1)^2 (1/{256}).
\end{equation*}
Hence, if at the same time $A^{\prime}(x_1)(1/4)\geq 2\left( \sqrt{2M_1/3} +
2+\sqrt{M_2/K}\right)$ and $KA^{\prime}(x_1)^2(1/256)>M_2$, we would get the
contradiction $M_2\geq \int_IP(A,B)dx>M_2$. This shows that 
\begin{equation*}
A^{\prime}(x_1)\leq 8\left(\sqrt{2M_1/3} +2+2\sqrt{M_2/K}\right)
\end{equation*}
if $x_1\leq 0$ and $A^{\prime}(x_1)\geq 0$. More generally, for every closed
interval $I$ of length $1$, there exists $x\in I$ such that 
\begin{equation*}
  P(A(x),B(x))\leq M_2,\quad \|(A(x),B(x))\|\leq 1+\sqrt {M_2/K},
\end{equation*}
and
\begin{equation*}
  |A^{\prime}(x)|\leq 8\left(\sqrt{2M_1/3} +2+2\sqrt{M_2/K}\right),
\end{equation*}
with $M_1,M_2\in[0,\infty)$ satisfying $\int_{I}\frac12|A^{\prime\prime}(x)|^2dx\leq M_1$ and $\int_I P(A(x),B(x))dx\leq M_2$.

Given $I$ as above, and a constant $\kappa >0$, observe
that, for all $\mu >0$, there exists $\nu >0$ such that, for all $(A,B)\in
H^{2}(\overset{\circ }{I})\times H^{1}(\overset{\circ }{I})$, the inequalities 
\begin{equation}
\int_{I}\left( \frac{\epsilon }{2}|A^{\prime \prime }|^{2}+\frac{1}{2
}|B^{\prime }|^{2}\right) dx<\kappa ~~\text{ and }~~\int_{I}P(A,B)dx<\nu
\label{eq: mu-nu}
\end{equation}
imply that
\begin{equation*}
\max_{I}\Big(P(A,B)+|A^{\prime }|\Big)<\mu .
\end{equation*}
This will be proven ad absurdum by assuming the opposite. There would exist $
\mu >0$ such that, for all integers $n\geq 1$, one could find $
(A_{n},B_{n})\in H^{2}(\overset{\circ }{I})\times H^{1}(\overset{\circ }{I})$
such that 
\begin{equation*}
\int_{I}P(A_{n},B_{n})dx<1/n~~\text{ and }~~\max_{I}\Big(
P(A_{n},B_{n})+|A_{n}^{\prime }|\Big)\geq \mu .
\end{equation*}
From this, one also gets 
\begin{equation*}
\min_{I}\|(A_{n},B_{n})\|\leq 1+\sqrt{1/(Kn)}~~\text{ and }
~~\min_{I}|A_{n}^{\prime }|\leq 8\left( \sqrt{2\kappa /(3\epsilon )}+2+2
\sqrt{1/(Kn)}\right) .
\end{equation*}
Hence the sequence $\{(A_{n},B_{n})\}$ is bounded in $H^{2}(\overset{\circ }{
I})\times H^{1}(\overset{\circ }{I})$. Taking a subsequence instead if
needed, it converges weakly in $H^{2}(\overset{\circ }{I})\times H^{1}(
\overset{\circ }{I})$, and thus strongly in $C^{1}(I)\times C(I)$, to some $
(A,B)\in H^{2}(\overset{\circ }{I})\times H^{1}(\overset{\circ }{I})$. One
gets the contradiction 
\begin{equation*}
\int_{I}P(A,B)dx=0~~\text{ and }~~\max_{I}\Big(P(A,B)+|A^{\prime }|\Big)\geq
\mu .
\end{equation*}
As a consequence, if $(A,B)\in H_{loc}^{2}(\mathbb{R})\times H_{loc}^{1}(
\mathbb{R})$ satisfies $J_{\epsilon }(A,B)<\infty $, then 
\begin{equation*}
\lim_{x\rightarrow \infty }P(A(x),B(x))=\lim_{x\rightarrow -\infty
}P(A(x),B(x))=0
\end{equation*}
and therefore the two limits $\lim_{x\rightarrow \infty }(A(x),B(x))$ and $
\lim_{x\rightarrow -\infty }(A(x),B(x))$ exist and belong to the set $\{(\pm 1,0),(0,\pm1)\}$. In the same way 
\begin{equation*}
\lim_{|x|\rightarrow \infty }A^{\prime }(x)=0.
\end{equation*}

\subsection{Minimizing sequences}

Let $\{(A_{n},B_{n})\}\subset X$ be a minimizing
sequence of $J_{\epsilon }$. Taking a subsequence if needed, it can be
assumed to converge weakly in $H_{loc}^{2}(\mathbb{R})\times H_{loc}^{1}(\mathbb{R})$, and strongly in $C_{loc}^{1}(\mathbb{R})\times C_{loc}(\mathbb{R})$, to some $(A,B)\in H_{loc}^{2}(\mathbb{R})\times H_{loc}^{1}(\mathbb{R})$ such that
\begin{equation}
\int_{\mathbb{R}}\Big(\frac{\epsilon }{2}|A^{\prime \prime }|^{2}+\frac{1}{2}
|B^{\prime }|^{2}+P(A,B)\Big)dx\leq \inf_{X}J_{\epsilon }.
\label{eq: couple limite}
\end{equation}
As it is unknown yet whether \eqref{e:limits} is satisfied, it is not yet possible to replace the inequality in the above formula by an equality.

After possible translations in $x$, one can suppose that $B_{n}(0)=1/2$ for
all $n\in \mathbb{N}$, because $(A_{n},B_{n})\in X$, and thus $B(0)=1/2$. It
remains to show that, up to a subsequence, 
\begin{equation*}
\lim_{x\rightarrow -\infty }(A(x),B(x))=(1,0)~~\text{ and }
~~\lim_{x\rightarrow \infty }(A(x),B(x))=(0,1).
\end{equation*}
Observe that 
\begin{equation*}
\lim_{n\rightarrow \infty}\int_{0}^{1}P(A_{n},B_{n})dx=\int_{0}^{1}P(A,B)dx>0.
\end{equation*}
We use a compactness by concentration result (see the appendix in \cite{BuGrWa} inspired by \cite{Be-Ce}).

Consider the Hilbert space $H=L^{2}((0,1))$. Define $\{u_{n}\}_{n \geq
1}\subset l^2(\mathbb{Z},H)$ by $u_n=(u_{n,j})_{j\in \mathbb{Z}}$ with 
\begin{equation*}
u_{n,j}=P^{1/2}(A_n(\cdot+j),B_n(\cdot+j)).
\end{equation*}
For $w\in \mathbb{Z}$, one denotes by $T_w:l^2(\mathbb{Z},H)\rightarrow l^2(\mathbb{Z},H)$ the translation operator $T_w (v_j)=(v_{j-w})$. Then the sequence $\{u_n\}$ is bounded in $l^2(\mathbb{Z},H)$, the set $\{u_{n,j}:j\in \mathbb{Z},n\geq 1\}$ is relatively compact in $H$ and 
\begin{equation*}
\liminf_{n\rightarrow \infty}\|u_{n,0}\|_{H}>0.
\end{equation*}

Given $\delta>0$, there exist an integer $k\geq 1$, $u^1,\ldots,u^k\in l^2(\mathbb{Z},H)$, and sequences $\{w^{\ell}_n\}_{n\in\mathbb{N}}\subset \mathbb{Z}$, $\ell=1,\ldots k$, such that 
\begin{equation*}
u^1\neq 0,\ldots,u^{k}\neq 0,
\end{equation*}
\begin{equation}  \label{eq: less epsilon}
\limsup_{n\rightarrow \infty}\left\| u_n-\sum_{\ell=1}^kT_{w^\ell_n}u^\ell
\right\|_{l^\infty(\mathbb{Z},H)}\leq \delta,
\end{equation}
\begin{equation*}  \label{eq: if k=1}
\lim_{n\rightarrow \infty}\left\|u_n-T_{w^1_n}u^1 \right\|_{l^\infty(\mathbb{
Z},H)}=0~~\text{ if }~~k=1,
\end{equation*}
after taking a subsequence of $\{u_n\}$ if needed. Furthermore,
\begin{equation}  \label{eq: distances}
\lim_{n\rightarrow \infty}|w^{\ell}_n-w^{\ell^{\prime}}_n|=\infty~~\text{
for all }~~ 1\leq \ell<\ell^{\prime}\leq k,
\end{equation}
if $k\geq 2$, 
\begin{equation}  \label{eq: weak limit for k'}
T_{-w^{k^{\prime}}_n}u_n\rightharpoonup u^{k^{\prime}},
\end{equation}
weakly in $l^2(\mathbb{Z},H)$, for each $1\leq k^{\prime}\leq k$
and 
\begin{equation}  \label{eq: inequality for sum}
\sum_{\ell=1}^k\|u^\ell\|_{l^2(\mathbb{Z},H)}^2 \leq\lim_{n\rightarrow
\infty}\|u_n\|_{l^2(\mathbb{Z},H)}^2,
\end{equation}
where the limit exists up to a subsequence. If $k\geq 2$, taking a subsequence
if needed and relabelling $u^1,\ldots,u^k$, one can also assume that 
\begin{equation*}
w^1_n<\ldots<w^k_n,\quad\forall n\in \mathbb{N}.
\end{equation*}

By \eqref{eq: weak limit for k'}, up to a subsequence, there exists, for $\ell\in\{1,\ldots,k\}$, $(A^\ell,B^\ell)\in H^2_{loc}(\mathbb{R})\times H^1_{loc}(\mathbb{R})$ such that 
\begin{equation}  \label{eq: quasi-def B^l}
(A_n(\cdot+w^{\ell}_n),B_n(\cdot+w^{\ell}_n)) \rightarrow (A^{\ell},B^{\ell})
\end{equation}
weakly in $H^2_{loc}(\mathbb{R})\times H^1_{loc}(\mathbb{R})$ and strongly
in $C^1_{loc}(\mathbb{R})\times C_{loc}(\mathbb{R})$, and 
\begin{equation*}
u^\ell_j=P^{1/2}(A^{\ell}(\cdot+j),B^\ell(\cdot+j)),\quad\forall j\in \mathbb{Z},
\end{equation*}
for $\ell\in\{1,\ldots,k\}$. 
The equation \eqref{eq: inequality for sum} gives 
\begin{equation*}
\sum_{\ell=1}^k\int_{\mathbb{R}}P(A^\ell,B^\ell)dx \leq \lim_{n\rightarrow
\infty}\int_{\mathbb{R}}P(A_n,B_n)dx.
\end{equation*}
Moreover, for all $1\leq k^{\prime}\leq k$, 
\begin{multline*}
\int_{\mathbb{R}}(B^{k^{\prime}})^{\prime}\left(B^{\prime}_n(\cdot+w^{k^{
\prime}}_n)-\sum_{\ell=1}^{k^{\prime}}(B^\ell)^{\prime}(\cdot+w^{k^{
\prime}}_n-w^\ell_n)\right)dx \\
=\int_{\mathbb{R}}(B^{k^{\prime}})^{\prime}B^{\prime}_n(\cdot+w^{k^{
\prime}}_n)dx -\int_{\mathbb{R}}(B^{k^{\prime}})^{\prime}(B^{k^{\prime}})^{
\prime}dx -\sum_{\ell=1}^{k^{\prime}-1}\int_{\mathbb{R}}(B^{k^{\prime}})^{
\prime}(B^\ell)^{\prime}(\cdot+w^{k^{\prime}}_n-w^\ell_n)dx \rightarrow 0,
\end{multline*}
by \eqref{eq: distances} and \eqref{eq: quasi-def B^l}, which implies that 
\begin{multline*}
\hspace*{-1ex}\lim_{n\rightarrow
\infty}\left\|B^{\prime}_n-\sum_{\ell=1}^{k^{\prime}-1}(B^\ell)^{\prime}(
\cdot-w^\ell_n)\right\|_{L^2(\mathbb{R})}^2 =\lim_{n\rightarrow
\infty}\left\|(B^{k^{\prime}})^{\prime}(\cdot-w^{k^{\prime}}_n)+\left(B^{
\prime}_n-\sum_{\ell=1}^{k^{\prime}}(B^\ell)^{\prime}(\cdot-w^\ell_n)\right)
\right\|_{L^2(\mathbb{R})}^2 \\
=\lim_{n\rightarrow
\infty}\left\|(B^{k^{\prime}})^{\prime}+\left(B^{\prime}_n(\cdot+w^{k^{
\prime}}_n)-\sum_{\ell=1}^{k^{\prime}}(B^\ell)^{\prime}(\cdot+w^{k^{
\prime}}_n-w^\ell_n)\right)\right\|_{L^2(\mathbb{R})}^2 \\
=\left\|(B^{k^{\prime}})^{\prime}\right\|_{L^2(\mathbb{R})}^2
+\lim_{n\rightarrow
\infty}\left\|B^{\prime}_n-\sum_{\ell=1}^{k^{\prime}}(B^\ell)^{\prime}(
\cdot-w^\ell_n)\right\|_{L^2(\mathbb{R})}^2
\end{multline*}
and thus 
\begin{equation*}
\sum_{k^{\prime}=1}^k\left\|(B^{k^{\prime}})^{\prime}\right\|_{L^2(\mathbb{R}
)}^2 \leq\lim_{n\rightarrow \infty}\left\|B^{\prime}_n\right\|_{L^2(\mathbb{R
})}^2\,.
\end{equation*}
In the same way 
\begin{equation*}
\sum_{k^{\prime}=1}^k\left\|(A^{k^{\prime}})^{\prime}\right\|_{L^2(\mathbb{R}
)}^2 \leq\lim_{n\rightarrow \infty}\left\|A^{\prime}_n\right\|_{L^2(\mathbb{R
})}^2
\end{equation*}
and 
\begin{equation*}
\sum_{k^{\prime}=1}^k\left\|(A^{k^{\prime}})^{\prime\prime}\right\|_{L^2(
\mathbb{R})}^2 \leq\lim_{n\rightarrow
\infty}\left\|A^{\prime\prime}_n\right\|_{L^2(\mathbb{R})}^2\,.
\end{equation*}
Hence 
\begin{equation}  \label{eq: sous-additif}
\sum_{\ell=1}^k J_\epsilon(A^{\ell},B^{\ell})\leq
\lim_{n\rightarrow\infty}J_\epsilon(A_n,B_n) = \inf_{X}J_\epsilon\,.
\end{equation}
From \eqref{eq: less epsilon} and the fact that $\lim_{|j|\rightarrow
\infty} \|u^\ell_j\|_H=0$ for all $\ell\in\{1,\ldots,k\}$, one gets 
\begin{equation*}
\sup\,\left\{\left\|u_{n,j}\right\|_{H} :\,j\in \mathbb{Z},\, |j-w^1_n|>
p,\,\ldots\,,|j-w^k_n| > p\right\}
\end{equation*}
\begin{equation*}
\leq \left\|u_n-\sum_{\ell=1}^kT_{w^\ell_n}u^\ell \right\|_{l^\infty(\mathbb{
Z},H)} +\sup\,\left\{\sum_{\ell=1}^k\left\|u^{\ell}_{j}\right\|_{H} :\,j\in 
\mathbb{Z},\, |j|> p \right\}
\end{equation*}
for each $n$ and 
\begin{equation*}
\lim_{p\rightarrow \infty}\, \limsup_{n\rightarrow \infty}\,
\sup\,\left\{\left\|u_{n,j}\right\|_{H} :\,j\in \mathbb{Z},\, |j-w^1_n|>
p,\,\ldots\,,|j-w^k_n| > p\right\}\leq \delta.
\end{equation*}
Given $\mu>0$ and $\epsilon$, one chooses $\kappa=2\inf_X J_\epsilon$ in 
\eqref{eq: mu-nu} and then $\delta=\sqrt{\nu/2}$, with $\nu>0$ as in \eqref{eq: mu-nu}, which gives 
\begin{equation*}
\lim_{p\rightarrow \infty}\, \limsup_{n\rightarrow \infty}\, \sup\,\left\{
P(A_n(x),B_n(x)),\; x\not\in [w^j_n-p,w^j_n+p+1],\; j=1,\dots,k \right\}\leq \mu.
\end{equation*}

Let $\rho>0$ be such that the open set $\{(a,b)\in \mathbb{R}^2:P(a,b)<\rho\}$ is the union of four open sets $V_{(0,\pm1)}$ and $V_{(\pm1,0)}$ with disjoint adherence, containing the points $(0,\pm1)$ and $(\pm1,0)$, respectively.
One can also suppose that the line $\mathbb{R}\times\{1/2\}$ does not meet $\{(a,b)\in \mathbb{R}^2:P(a,b)\leq \rho\}$.

If one chooses $\mu=\rho/2$, then $p$ large enough, one gets for all $n$
large enough, up to a subsequence, 
\begin{equation*}
P(A_n(x),B_n(x))<\rho,\quad \forall x \in(-\infty, w^1_n-p),
\end{equation*}
\begin{equation*}
P(A_n(x),B_n(x))<\rho,\quad\forall x \in(w^k_n+p+1,\infty),
\end{equation*}
and, if $k\geq 2$, 
\begin{equation*}
P(A_n(x),B_n(x))<\rho,\quad\forall x \in(w^{\ell_n-1}+p+1,w^\ell_n-p)\neq\emptyset,
\end{equation*}
for all $\ell\in\{2,\ldots,k\}$. In the two first cases, as well as in the
last case for each $\ell\in\{2,\ldots,k\}$, $(A_n(x),B_n(x))$ not only
satisfies $P(A_n(x),B_n(x))<\rho$, but $(A_n(x),B_n(x))$ even stays in $
V_{(0,1)}$, $V_{(0,-1)}$, $V_{(1,0)}$ or $V_{(-1,0)}$ (this can change in
each case and one uses the continuity of $(A_n,B_n)$). As $(A_n,B_n)\in X$,
one has the following additional informations: 
\begin{equation*}
(A_n(x),B_n(x))\in V_{(1,0)},\quad\forall x \in(-\infty, w^1_n-p),
\end{equation*}
and 
\begin{equation*}
(A_n(x),B_n(x))\in V_{(0,1)},\quad\forall x \in(w^k_n+p+1,\infty).
\end{equation*}
Hence there exists $\hat \ell\in\{1,\ldots,k\}$ such that 
\begin{equation*}
(A_n(x),B_n(x))\in V_{(1,0)}\cup V_{(-1,0)},\quad\forall x \in(w^{\hat \ell-1}_n+p+1, w^{\hat \ell}_n-p),
\end{equation*}
and 
\begin{equation*}
(A_n(x),B_n(x))\in V_{(0,1)}\cup V_{(0,-1)},\quad\forall x \in(w^{\hat \ell}_n+p+1, w^{\hat \ell+1}_n-p),
\end{equation*}
with the understanding that $w^{0}_n=-\infty$ and $w^{k+1}_n=+\infty$. From 
\eqref{eq: weak limit for k'}, it follows that 
\begin{equation*}
\lim_{x\rightarrow -\infty}
(A^{\hat\ell}(x),B^{\hat\ell}(x))\in\{(\pm1,0)\} ~~\text{ and }~~
\lim_{x\rightarrow \infty}
(A^{\hat\ell}(x),B^{\hat\ell}(x))\in\{(0,\pm1)\}.
\end{equation*}
As, with the right choice of signs, one has $(\pm A^{\hat\ell},\pm B^{\hat
\ell})\in X$ and 
\begin{equation*}
\inf_{X}J_\epsilon\leq J_\epsilon(\pm A^{\hat\ell},\pm B^{\hat
\ell})=J_\epsilon(A^{\hat\ell},B^{\hat \ell}),
\end{equation*}
one gets (see \eqref{eq: sous-additif}) 
\begin{equation*}
\inf_{X}J_\epsilon\leq \sum_{\ell=1}^k J_\epsilon(A^{\ell},B^{\ell})\leq
\lim_{n\rightarrow\infty}J_\epsilon(A_n,B_n) = \inf_{X}J_\epsilon\,.
\end{equation*}
As $u^\ell\neq 0$ for all $\ell\in\{1,\ldots,k\}$, this is only possible if $k=1$ and 
\begin{equation*}
J_\epsilon(\pm A^{\hat\ell},\pm B^{\hat
\ell})=J_\epsilon(A^{\hat\ell},B^{\hat \ell}) =\inf_{X}J_\epsilon.
\end{equation*}
Since $k=\hat\ell =1$, one also has 
\begin{equation*}
\lim_{x\rightarrow -\infty} (A^{1}(x),B^{1}(x))=(1,0) ~~\text{ and }~~
\lim_{x\rightarrow \infty} (A^{1}(x),B^{1}(x))=(0,1)
\end{equation*}
This shows that $(A^{1},B^{1})\in X$ minimizes $J_\epsilon$. In addition, up
to a translation in $x$, it is equal to $(A,B)\in H^2_{loc}(\mathbb{R})\times H^1_{loc}(\mathbb{R})$ introduced in \eqref{eq: couple limite}, which therefore indeed belongs to $X$.

Finally, notice that $(A,|B|)$ is also a minimal pair and therefore one can assume that $B\geq 0$  on $\mathbb{R}$. 
However $(A,B)$ tends to $(0,1)$ as $x\rightarrow \infty $ in a way such that $A$ oscillates around $0$. This behavior is given by the linearization at $(A,B)=(0,1)$, because this equilibrium is hyperbolic.

\subsection{Limit $g\to 1$}

It remains to prove the last property in Theorem~\ref{t:main}. We show below that any heteroclinic orbit connecting the equilibria $(1,0)$ and $(0,1)$ which minimizes the functional $J_\epsilon$ in the space $X$ remains in a neighborhood of the circle $A^2+B^2=1$ as $g\rightarrow 1+$.

Let us first estimate $\min_X J_\epsilon$ as $g\rightarrow 1+$.
We introduce the ``test function''
$(A_1,B_1)$ defined as follows:
$$(A_1,B_1)(x)=\left(\cos\left(\frac \pi 4+\frac{\arctan(x)}2\right),\ \sin\left(\frac \pi 4+\frac{\arctan(x)} 2\right)\right),\quad x\in \R.$$
Then $(A_1,B_1)$ belongs to the space $X$, and we have the formulas for the first and second order derivatives:
$$(A'_1,B'_1)(x)=\left(-\sin\left(\frac \pi 4+\frac{\arctan(x)}2\right),\cos\left(\frac \pi 4+\frac{\arctan(x)} 2\right)\right)\frac 1 {2(x^2+1)}\,,$$
$$(A''_1,B''_1)(x)=\left(-\cos\left(\frac \pi 4+\frac{\arctan(x)}2\right),-\sin\left(\frac \pi 4+\frac{\arctan(x)} 2\right)\right)\frac 1 {4(x^2+1)^2}
$$
$$
-\left(-\sin\left(\frac \pi 4+\frac{\arctan(x)}2\right),\cos\left(\frac \pi 4+\frac{\arctan(x)} 2\right)\right)\frac x {(x^2+1)^2}.$$
For a suitably chosen positive constant $C$, we obtain the estimates:
\begin{eqnarray*}
  && | B_1(x)-1|\leq C\,\frac 1 x, \quad \forall\ x\leq -1,\\
  && 0<B_1(x)\leq C\,\frac 1 x,\quad \forall\ x\geq 1,\\
  &&0<A_1(x)\leq C\,\frac 1 x,\quad \forall\ x\leq -1,\\
  &&|A_1(x)-1|\leq C\,\frac 1 x,\quad \forall\ x\geq 1.
\end{eqnarray*}
For $\gamma>0$, let $(A_\gamma,B_\gamma)\in X$ be defined by
$$(A_\gamma,B_\gamma)(x)= (A_1(\gamma x),B_1(\gamma x)),\quad\forall\ x\in \R.$$
Then
$$\min_{X} J_{\epsilon,g}\leq J_{\epsilon,g} (A_\gamma,B_\gamma)
= \int_{\mathbb R}\Big(\frac \epsilon 2 |A_{\gamma}''|^2+\frac 1 2 |B_{\gamma}'|^2+\frac {g-1} 2 A_{\gamma}^2B_{\gamma}^2\Big)d x$$
$$= \int_{\mathbb R}\Big(\frac{\gamma^3 \epsilon} 2 |A_{1}''|^2+\frac {\gamma} 2 |B_{1}'|^2+\frac {\gamma^{-1}(g-1)} 2 A_{1}^2B_{1}^2\Big)d x.$$
By choosing $\gamma=(g-1)^{1/2}$, we get
$\min_{X} J_{\epsilon,g}\rightarrow 0$
as $g\rightarrow 1+$ (for fixed $\epsilon>0$).

Let now $(A_g,B_g)$ denote a minimizing heteroclinic orbit, where $\epsilon>0$ is still fixed, but  we insist on the dependence on the parameter $g>1$.
As we have seen, for some constant $K>0$ and each closed interval $I$ of length $1$, there exists $x_0\in I$ such that
$$
P(A_g(x_0),B_g(x_0))\leq\min_X J_{\epsilon,g},\quad
\|(A_g(x_0),B_g(x_0))\|\leq 1+\sqrt {\min_X J_{\epsilon,g}/K}$$
and
$$|A_g'(x_0)|\leq 
8\left(\sqrt{2\min_X J_{\epsilon,g}/(3\epsilon)} +2+2\sqrt{\min_X J_{\epsilon,g}/ K}\right).$$

For $\tilde \gamma>0$, consider an non-empty open interval $\tilde I\subset I$, if any, such that
\[
\frac 1 4 (A_g(x)^2+B_g(x)^2-1)^2\geq \tilde\gamma >0,\quad \forall\ x\in\tilde I.
\]
Then its length  $|\tilde I|$ satisfies $|\tilde I|\leq \tilde\gamma^{-1}\min_X J_{\epsilon,g}$ because
$$J_{\epsilon,g}(A_g,B_g)\geq \int_{\tilde I}P(A_g,B_g)dx\geq |\tilde I|\tilde \gamma.$$
Moreover, for all $x_1,x_2,x_3\in \tilde I$, 
$$|B_g(x_1)-B_g(x_2)|\leq \int_{\tilde I}|B_g'(x)|dx
\leq |\tilde I|^{1/2} (2\min_X J_{\epsilon,g})^{1/2}
\leq \tilde\gamma^{-1/2}2^{1/2}\min_X J_{\epsilon,g}\,,
$$
$$|A_g'(x_3)-A_g'(x_0)|\leq \int_{I}|A_g''(x)|dx
\leq |I|^{1/2} (2\epsilon^{-1}\min_X J_{\epsilon,g})^{1/2}
= (2\epsilon^{-1}\min_X J_{\epsilon,g})^{1/2},
$$
$$|A_g'(x_3)|\leq 
8\left(\sqrt{2\min_X J_{\epsilon,g}/(3\epsilon)} +2+2\sqrt{\min_X J_{\epsilon,g}/K}\right)
+ (2\epsilon^{-1}\min_X J_{\epsilon,g})^{1/2}$$
and
\begin{eqnarray*}
|A_g(x_1)-A_g(x_2)|&\leq& 
8|\tilde I|
\left(\sqrt{2\min_X J_{\epsilon,g}/(3\epsilon)} +2+2\sqrt{\min_X J_{\epsilon,g}/K}\right)\\
&&+|\tilde I| (2\epsilon^{-1}\min_X J_{\epsilon,g})^{1/2}=O(\tilde\gamma^{-1}\min_X J_{\epsilon,g}).
\end{eqnarray*}
Setting $\tilde\gamma=(\min_X J_{\epsilon,g})^{1/2}$, we get $|\tilde I|<1$ for $g-1>0$ small enough. Moreover, for all $x_2\in I$, there exists $x_1\in I$ such that
$$\frac 1 4 (A_g(x_1)^2+B_g(x_1)^2-1)^2\leq\tilde\gamma=(\min_X J_{\epsilon,g})^{1/2}\rightarrow 0$$
and
$$\|(A_g,B_g)(x_2)-(A_g,B_g)(x_1)\|=O((\min_X J_{\epsilon,g})^{1/2})\rightarrow 0,$$
as $g\rightarrow 1 +$.
Indeed, if 
$\frac 1 4 (A_g(x_2)^2+B_g(x_2)^2-1)^2\leq \tilde \gamma$, choose $x_1=x_2$ and, if
$\frac 1 4 (A_g(x_2)^2+B_g(x_2)^2-1)^2 >\tilde \gamma$, choose for $x_1$ one of the boundary points of the largest open interval $\tilde I$ included in $I$, containing $x_2$ and on which $\frac 1 4 (A_g(x)^2+B_g(x)^2-1)^2> \tilde\gamma$.

Consequently,
\begin{eqnarray*}
  \frac 1 4 (A_g(x_2)^2+B_g(x_2)^2-1)^2&=&\frac 1 4 \left(A_g(x_1)^2+B_g(x_1)^2-1+O((\min_X J_{\epsilon,g})^{1/2})\right)^2\\
  &=&O((\min_X J_{\epsilon,g})^{1/2})
\end{eqnarray*}
and, the estimates being uniform with respect to any arbitrary interval $I$ of length $1$, we conclude that
$$\lim_{g\rightarrow 1+}\sup_{x\in \R}\frac 1 4 (A_g(x)^2+B_g(x)^2-1)^2=0.$$
This completes the proof of Theorem~\ref{t:main}.

\end{document}